\documentclass{elsart_mod}

\usepackage{amssymb}
\usepackage{amsmath} 
\usepackage{mathrsfs}

\newenvironment{nonnumberedtheorem}[1]{\begin{trivlist}\item[]{\bf #1\/} 
  \em\ignorespaces }{\end{trivlist}}

\begin{document}

\newcommand{\aut}[1]{\mbox{\rm Aut}(#1)}
\newcommand{\heins}[2]{\mbox{\rm H}^{1}(#1,#2)}
\newcommand{\UNIT}[1]{\mbox{\rm U}(#1)}
\newcommand{\cen}[2]{\mbox{\rm C}_{#1}(#2)}
\newcommand{\nor}[2]{\mbox{\rm N}_{#1}(#2)}
\newcommand{\zen}[1]{\mbox{\rm Z}(#1)}
\newcommand{\cyk}[1]{C_{#1}}
\newcommand{\opn}[2]{\mbox{\rm O}_{#1}(#2)}
\newcommand{\opns}[1]{\mbox{\rm O}_{p^{\prime}}(#1)}
\newcommand{\lara}[1]{\langle{#1}\rangle}
\newcommand{\ZZ}{\mbox{$\mathbb{Z}$}}
\newcommand{\sZZ}{\mbox{$\scriptstyle\mathbb{Z}$}}
\newcommand{\QQ}{\mbox{$\mathbb{Q}$}}
\newcommand{\NN}{\mbox{$\mathbb{N}$}}
\newcommand{\AUG}[2]{\mbox{\rm I}_{#1}(#2)}
\newcommand{\NU}[1]{\mbox{\rm V}(#1)}
\newcommand{\ol}[1]{\overline{#1}}

\begin{frontmatter}


\title{Units of $\boldsymbol{p}$-power order in principal 
$\boldsymbol{p}$-blocks of $\boldsymbol{p}$-constrained groups}
\runtitle{Units of $p$-power order in principal $p$-blocks}
\thanks[]{Research supported by the Deutsche Forschungsgemeinschaft.}

\author{Martin Hertweck\thanksref{}}
\ead{hertweck@mathematik.uni-stuttgart.de}
\address{Universit\"at Stuttgart, Fachbereich Mathematik,
IGT, 70550 Stuttgart, Germany}

\begin{abstract}
Let $G$ be a finite group having a normal $p$-subgroup $N$ that
contains its centralizer $\cen{G}{N}$, and let $R$ be a $p$-adic
ring. It is shown that any finite $p$-group of units of augmentation one
in $RG$ which normalizes $N$ is conjugate to a subgroup of $G$ by a 
unit of $RG$, and if it centralizes $N$ it is even contained in $N$.
\end{abstract}

\end{frontmatter}


\section{Introduction}

This paper grew out of an attempt to understand more fully 
part of a theorem due to Roggenkamp and Scott (see
\cite[Theorem~6]{rogg:86}, \cite{sco:87,sco:90},
\cite[Theorem~19]{rogg:91b}) about conjugacy of certain
finite $p$-subgroups in the group of units of a $p$-adic group ring.
The theorem in question, 
by now called the F$^{\ast}$-Theorem, is stated below together 
with references where a detailed account on its proof can be found
(with the result from Section~2 of the present paper being relevant).

Some of the interesting aspects of the group of units of a
group ring $SG$ of a finite group $G$ concern its finite subgroups,
in particular when the coefficient ring $S$ is a $G$-adapted ring, i.e., 
an integral domain of characteristic zero in which no prime divisor of 
the order of $G$ is invertible.
Below, a few well known results in this case are listed.
Note that it suffices to consider only the group of units
$\NU{SG}$ consisting of units of augmentation one.
For $u\in\NU{SG}$, we say that $u$ is a trivial unit if $u\in G$,
and the trace of $u$ is its $1$-coefficient (with respect to the 
basis $G$). Let $S$ be a $G$-adapted ring. 
Then (see \cite{sak:71}, \cite{kar:89} or \cite{seh:93}):
\begin{enumerate}
\item[(a)]
a non-trivial unit of $SG$ of finite order has trace zero;
\item[(b)]
the order of a finite subgroup of $\NU{SG}$ divides the order of $G$;
\item[(c)]
a central unit of finite order is a trivial unit.
\end{enumerate}

Limiting attention only to finite $p$-subgroups in the group of units,
one might ask whether comparable results hold with $S$ replaced by the 
ring $\ZZ_{p}$ of $p$-adic integers. However, (a) and (b) does not 
carry over, even not if $\ZZ_{p}G$ consists of a single block only
(cf.\ \cite[Section~XIV]{RoTa:92}).
Imposing additional conditions one might also ask how certain
finite $p$-subgroups are embedded in $\NU{\ZZ_{p}G}$.
If, for example, attention is directed to the principal block $B$ of 
$\ZZ_{p}G$, a Sylow $p$-subgroup $P$ of $G$ is identified with its
projection on $B$, and $\alpha$ is an augmented automorphism of $B$,
then the question whether $P$ is conjugate by a unit of $B$ 
to its image $P\alpha$, is part of Scott's 
``defect group (conjugacy) question'' 
(see \cite[p.~267]{sco:87}, \cite{sco:90}).
For $p$-groups $G$, this question was answered in the affirmative
by Roggenkamp and Scott \cite{RoSc:87a}.

Here, the following two theorems are proved.
In both we assume that $G$ has a normal $p$-subgroup $N$ satisfying 
$\cen{G}{N}\leq N$.
By definition, this means that $G$ is $p$-constrained and $\opns{G}=1$ 
(see \cite[VII,~13.3]{HuBl:82a}).

Throughout the paper, $R$ denotes a $p$-adic ring, that is, the
integral closure of the $p$-adic integers $\ZZ_{p}$ in a finite extension
field of the $p$-adic field $\QQ_{p}$. 
(Then $R$ is a complete discrete valuation ring.)
Note that by our assumption on $G$, the group ring $RG$ will consist
of a single (principal) block only (see \cite[VII,~13.5]{HuBl:82a}). 

\begin{nonnumberedtheorem}{Theorem~A.}
Suppose that $G$ has a normal $p$-subgroup $N$ that contains its
centralizer $\cen{G}{N}$.
Then any finite $p$-group in $\NU{RG}$ which normalizes $N$ 
is conjugate to a subgroup of $G$ by a unit of $RG$.
\end{nonnumberedtheorem}

\begin{nonnumberedtheorem}{Theorem~B.}
Suppose that $G$ has a normal $p$-subgroup $N$ that contains its
centralizer $\cen{G}{N}$. Then any finite $p$-group in $\NU{RG}$ 
which centralizes $N$ is contained in $N$.
\end{nonnumberedtheorem}

The proofs are somewhat complicated by the fact that we do not know in 
advance that $RG$ is free for the ``multiplication action'' of the finite 
$p$-group under consideration (taking this for granted, Theorem~A should
be part of the F$^{\ast}$-Theorem).
Section~2 contains some preparatory results needed for the handling
of the case $p=2$. Theorems~A and B are proved in Section~3.
The bimodule arguments used there are inspired by \cite[p.~231]{rogg:86}.
The proof depends heavily on the strong results of Weiss
on $p$-permutation lattices (see \cite{wei:87,wei:93}, \cite{rogg:92b}).
That these results can be applied rests upon 
the ``Ward--Coleman Lemma.''
Coleman's contribution \cite{col:64} is well known, but the first version
of the lemma appears in an article of Ward \cite{War:60} as a contribution
to a seminar run by Richard Brauer at Harvard. 
See also \cite[Proposition~1.14]{sak:71}, \cite[2.6~Theorem]{JaMa:87}.
Actually, for its proof it is only needed that $p$ is not invertible 
in the commutative ring $R$.

\begin{nonnumberedtheorem}{Ward--Coleman Lemma.}\label{gcr}
Let $H$ be a $p$-subgroup of the finite group $G$. Then
$\mbox{\rm N}_{\text{\rm V}(RG)}(H)=\mbox{\rm N}_{G}(H)\cdot
\mbox{\rm C}_{\text{\rm V}(RG)}(H)$.
\qed
\end{nonnumberedtheorem}

A consequence of Theorem~B deserves explicit mention.
\begin{nonnumberedtheorem}{Corollary.}
Suppose that $G$ has a normal $p$-subgroup $N$ that contains its
centralizer $\cen{G}{N}$.
Then each central unit of (finite) $p$-power order in $\NU{RG}$ 
is a trivial unit, i.e.\ contained in $\zen{G}$.
\end{nonnumberedtheorem}

It is not known whether a corresponding result holds for the principal
block of an arbitrary (non-solvable) group. Progress in this direction
might lead to applications in finite group theory,
through Robinson's work \cite{rob:90}
on odd analogues of Glauberman's Z$^{\ast}$-Theorem.

Finally, we state the theorem of Roggenkamp and Scott. 

\begin{nonnumberedtheorem}{F$^{\boldsymbol{\ast}}$-Theorem.}
Let $G$ be a finite group having a normal $p$-subgroup $N$ that
contains its centralizer $\cen{G}{N}$.
If $\alpha$ is an automorphism of $RG$ stabilizing both the
augmentation ideal $\AUG{R}{G}$ and the ideal $\AUG{R}{N}G$, then the
groups $G$ and $G^{\alpha}$ are conjugate by a unit of $RG$.
\end{nonnumberedtheorem}

The first step of the proof consists in showing that $N$ and 
$N^{\alpha}$ are conjugate by a unit of $RG$, so that $N=N^{\alpha}$ 
can be assumed (details of sketch proof in \cite{rogg:86}
are given in \cite[Lemma~4.1]{HeKi:00b}).
Then the image $P^{\alpha}$ of a Sylow $p$-subgroup $P$ of $G$
normalizes $N$, and by Theorem~A one can even assume that $P=P^{\alpha}$.
Now it can be shown that there exists an automorphism $\rho$ of $G$
such that the composition $\alpha\rho$ fixes the elements of $P$, 
which implies that $\alpha\rho$ is an inner automorphism of $RG$
(see \cite{He:04} for a published account).

\section{(Mod $\boldsymbol{2}$) conjugacy for involutions}

The main result of this section is Proposition~\ref{cp}, which will
be used in conjunction with Proposition~\ref{sbor} to prove in the 
next section Lemma~\ref{l2} for $p=2$, which then allows application of
Weiss' theorem on $p$-permutation lattices 
in the proof of Theorem~A and B.

We continue to denote by $G$ a finite group and by $R$ a $p$-adic ring.
For a subgroup $H$ of $G$, let $\AUG{R}{H}$ denote the 
augmentation ideal of $RH$, i.e.\ the $R$-span of the elements $h-1$
($h\in H$). If $H$ is a normal subgroup, $\AUG{R}{H}G$ is the kernel of 
the natural map $RG\rightarrow RG/H$.
For any subset $T$ of $G$, let $R[T]$ denote
the $R$-span of the elements of $T$.

The Ward--Coleman Lemma suggests a closer investigation of
$\mbox{\rm C}_{\mbox{\rm\scriptsize V}(RG)}(N)$, for $G$ and $N$ as in
Theorems~A and B.

\begin{prop}\label{sbor}
Suppose that $G$ has a normal $p$-subgroup $N$ that contains its
centralizer $\cen{G}{N}$. Then the following holds.
\begin{enumerate}
\item[{\rm (a)}]
$\cen{RG}{N}\subseteq R+\AUG{R}{N}G+pRG$.
\item[{\rm (b)}]
Assume that $p=2$ and set $\bar{G}=G/\opn{2}{G}$.
If a subset $T$ of $G$ is chosen such that $\bar{T}$ is the set 
of involutions in $\bar{G}$, then 
\[ \cen{RG}{N}\subseteq R+\AUG{R}{\opn{2}{G}}G+4RG+2R[T]. \]
\end{enumerate}
\end{prop}
\begin{pf}
For $g\in G$ define its ``$N$-class sum'' $\mathscr{C}_{g}^{N}$
to be the sum (in $RG$) of the distinct $g^{x}$, $x\in N$.
Note that the $N$-class sums form an $R$-basis of $\cen{RG}{N}$.
Let $g\in G$. If $\mbox{\rm C}_{N}(g)=N$, then $g\in\cen{G}{N}\leq N$ and 
$\mathscr{C}_{g}^{N}=g$, and if $|N:\mbox{\rm C}_{N}(g)|=p^{a}>1$, then
$\mathscr{C}_{g}^{N}\in p^{a}g + \AUG{R}{N}G$. Hence (a) holds.
Set $N^{\ast}=N/\Phi(N)$, the quotient of $N$ by its Frattini subgroup 
$\Phi(N)$. The group $G$ acts by conjugation on the elementary abelian 
group $N^{\ast}$. Since $N$ is not centralized by $p^{\prime}$-elements 
of $G\setminus\{1\}$, the kernel $K$ of this operation is a $p$-group 
(see \cite[III,~3.18]{hupp:67}), whence contained in $\opn{p}{G}$. 
Now let $p=2$ and suppose that $|N:\mbox{\rm C}_{N}(g)|=2$ for some
$g\in G$.
Then $|N^{\ast}:\cen{N^{\ast}}{g}|\leq 2$, so $g^{2}\in K\leq\opn{2}{G}$ and
$\mathscr{C}_{g}^{N}\in R+\AUG{R}{\opn{2}{G}}G+2R[T]$, which completes 
the proof of (b).
\qed
\end{pf}

In particular, we are led to consider units $u$ of order $p$ in $RG/N$ 
which map to $1$ in $(R/pR)(G/N)$. The following well known lemma 
(cf.\ \cite[Lemma~(53.3)]{wei:93}) shows that this is actually 
only a matter for $p=2$.

\begin{lem}\label{odd}
Let $p$ be an odd prime.
Then the group of units in $1+pRG$ has no $p$-torsion.
\end{lem}
\begin{pf}
Suppose that $u$ is a unit of order $p$ in $u\in 1+pRG$ and 
write $u=1+p^{a}x$ with $a\in\NN$ such that $x\in RG\setminus pRG$. 
Writing out the binomial expansion for $(1+p^{a}x)^{p}=1$ and solving 
for $x$, we get $p\cdot p^{a}x=-\sum_{n=2}^{p}\tbinom{p}{n}p^{an}x^{n}$.
But the right-hand side has a factor 
$p\cdot p^{2a}$ since $p>2$. This contradiction proves the lemma.
\qed
\end{pf}
The reader might wish to convince himself that the statement of the lemma
does not hold for $p=2$, even if $u$ is assumed to have augmentation one.

As usual, we let $[RG,RG]$ denote the 
additive commutator of $RG$, defined by
\[ [RG,RG]=\mbox{}_{R}\lara{g^{h}-g \mid g,h\in G}
=\{xy-yx \mid x,y\in RG\}. \]

The following formulas are well known (see, for example,
\cite[Lemma~(7.1)]{seh:93}).
We write $K$ for the quotient field of $R$, let
$a_{1},\ldots,a_{l}\in RG$ and $n\in\NN$.
\begin{align*}
& [KG,KG]\cap RG=[RG,RG] \tag{1.1}\label{1.1} \\
& [RG,RG]^{p}\subseteq pRG+[RG,RG] \tag{1.2}\label{1.2} \\
& (a_{1}+\cdots +a_{l})^{p^{n}}\equiv a_{1}^{p^{n}}+\cdots +a_{l}^{p^{n}}
\mbox{ {}\rm mod } pRG+[RG,RG] \tag{1.3}\label{1.3} \\
\intertext{Let $G_{p^{\prime}}$ denote the set of $p^{\prime}$-elements 
of $G$. It follows immediately from (1.2) and (1.3) that for 
large enough $n\in\NN$ (it suffices that $p^{n}$ is the order of a 
Sylow $p$-subgroup of $G$) and any $x\in RG$,}
& x^{p^{n}}\in pRG+[RG,RG]+R[G_{p^{\prime}}]. \tag{1.4}\label{1.4}
\end{align*}
\addtocounter{equation}{1}

Fix some $c\in G$. We define
\[ [RG]^{1-c}=\mbox{}_{R}\lara{g-g^{c}\mid g\in G}
=\{x-x^{c}\mid x\in RG\}. \]
Set $C=\lara{c}$.
If $RG$ is considered as a $C$-module for the ``conjugation action''
$m\cdot c=c^{-1}mc$ ($m\in RG$), then $[RG]^{1-c}$ is the image of the 
map $RG\xrightarrow{\cdot(1-c)} RG$.
Since $\heins{C}{RG}=0$ (see \cite[(1.4.1)~Proposition]{RoSc:87a}), 
$[RG]^{1-c}$ is also the kernel of the map
$RG\xrightarrow{\cdot\widehat{C}} RG$ where 
$\widehat{C}$ denotes the sum of the elements in $C$.
In particular, it follows that 
\begin{equation}
[KG]^{1-c}\cap RG=[RG]^{1-c}.
\end{equation}

We remark that $[RG]^{1-c}$ is invariant under left (and right) 
multiplication with $c$. 
Also, if $c^{2}=1$, then $([RG]^{c-1})^{2}\subseteq 2RG+[RG]^{c-1}$
in analogy to (1.2) holds since
$(g-g^{c})(h-h^{c})=gh+(gh)^{c}-g^{c}h-(g^{c}h)^{c}$ for all $g,h\in G$.
\begin{lem}\label{abc}
Let $c\in G$ and $u\in\NU{RG}$ be $2$-elements.
Assume that $c^{2}=u^{2}$ and $c^{-1}u=1+2f$ for some $f\in RG$. Then 
\begin{enumerate}
\item[{\rm (a)}]
$f+f^{2}\in 2RG+[RG]^{1-c}$;
\item[{\rm (b)}]
$f+f^{2^{n}}\in 2RG+[RG,RG]$ for all $n\in\NN$;
\item[{\rm (c)}]
$f\in 2RG+[RG,RG]+R[G_{2^{\prime}}]$.
\end{enumerate}
\end{lem}
\begin{pf}
Squaring the equation $c^{-1}u=1+2f$ yields
$1+4(f+f^{2})=(c^{-1}u)^{2} $, so $4(f+f^{2})u^{-1}c=c^{-1}u-u^{-1}c$.
Now $c^{-1}u=cu^{-1}$ since $c^{2}=u^{2}$, and it follows that
$4(f+f^{2})u^{-1}c=cu^{-1}-(cu^{-1})^{c}\in[RG]^{1-c}$.
By (2), and since $c^{-1}u\in 1+2RG$, (a) follows.
Now (b) follows inductively from (1.2) and
$f+f^{2^{n}}=(f+f^{2})+(f+f^{2^{n-1}})^{2}-2f(f+f^{2^{n-1}})$.
Part (c) follows from (b) and (1.4).
\qed
\end{pf}

Let $T_{2}$ be the set of involutions in $G$.
\begin{lem}\label{l1}
If $x\in R[T_{2}]\cap\AUG{R}{G}$ then 
$x^{2}\in (R[G\setminus T_{2}]\cap [RG,RG])+2RG$.
\end{lem}
\begin{pf}
Write $x=\sum_{s\in T_{2}}r_{s}(s-1)$ with $r_{s}\in R$ for all
$s\in T_{2}$. Then
\[ x^{2}=\sum_{s\in T_{2}}r_{s}^{2}(s-1)^{2}\;+\!\!
\sum_{\substack{\{s,t\}\subseteq T_{2} \\ s\neq t}}
r_{s}r_{t}((s-1)(t-1)+(t-1)(s-1)). \]
Let $s,t\in T_{2}$. Then $(s-1)^{2}=2(1-s)\in 2 RG$ and
\begin{align*}
& (s-1)(t-1)+(t-1)(s-1) \\
& \qquad =2(1-s-t+st)-st+ts \\
& \qquad \in
\begin{cases}
2RG & \text{if $[s,t]=1$,} \\
2RG+(R[G\setminus T_{2}]\cap [RG,RG]) & \text{otherwise.}
\end{cases}
\end{align*}
This proves the lemma.
\qed
\end{pf}

It is evident that up to now, letting $R$ be a $p$-adic ring
was unnecessarily restrictive, but finally we shall use results
of Weiss on permutation modules to prove:

\begin{prop}\label{cp}
Suppose that $p=2$, so that $R$ is a $2$-adic ring, and
let $c\in G$ and $u\in\NU{RG}$ be $2$-elements with
$c^{2}=u^{2}$ and $c^{-1}u=1+2f$ for some $f\in 2RG +R[ T_{2}] + R$.
Then $c$ and $u$ are conjugate in the units of $RG$.
\end{prop}
\begin{pf}
By Lemma~\ref{abc}(c), $f\in 2RG+[RG,RG]+R[G_{2^{\prime}}]$.
Compared with $f\in 2RG +R[ T_{2}] + R$, it follows that 
$f\in 2RG +(R[ T_{2}]\cap[RG,RG])+R$.
Taking augmentation gives
$f\in 2RG +(R[ T_{2}]\cap[RG,RG])$, so that by Lemma~\ref{l1},
$f^{2}\in 2RG+(R[G\setminus T_{2}]\cap [RG,RG])$.
By Lemma~\ref{abc}(a), $f+f^{2}\in 2RG+[RG]^{1-c}$. 
The last three equations show that 
\begin{equation}\label{eq1}
f\in 2RG+[RG]^{1-c}.
\end{equation}
Let $C$ be a cyclic group with generator $x$ of order the maximum 
of the orders of $c$ and $u$.
Consider $RG$ as an $RC$-lattice, the action of $x$
given by $m\cdot x=c^{-1}mu$ for all $m\in RG$. A $1$-cocycle 
$\beta : C\rightarrow RG$ is defined by $\beta(x)=f=(c^{-1}u-1)/2$. 
Notice that $c$ and $u$ are conjugate in the units of $RG$ if
$\beta$ is a $1$-coboundary: For if $\beta(x)=m-c^{-1}mu$ for some $m\in RG$,
then $v=1+2m$ is a unit in $RG$ (see \cite[(5.10)]{CuRe:81})
and $u=c^{v}$. The exact sequence 
$0\longrightarrow RG\stackrel{2}{\longrightarrow} RG
\stackrel{-}{\longrightarrow}\overline{RG}\longrightarrow 0$
has cohomology exact sequence
\[ \ldots \longrightarrow 
\heins{C}{RG} \stackrel{2}{\longrightarrow} \heins{C}{RG}
\stackrel{-}{\longrightarrow}\heins{C}{\overline{RG}}\longrightarrow\ldots \]
By \cite[Theorem~(50.2)]{wei:93}, $RG$ is a monomial lattice for $C$ over $R$,
so that $2$ annihilates $\heins{C}{RG}$ by \cite[Lemma~(53.1)]{wei:93}.
Hence the cohomology exact sequence implies that 
$\heins{C}{RG}\stackrel{-}{\longrightarrow}\heins{C}{\overline{RG}}$ is 
injective. By \eqref{eq1}, the class of $\beta$ maps to zero under this map; 
hence $\beta$ is a $1$-coboundary, and $c$ and $u$ are conjugate
in the units of $RG$. 
\qed
\end{pf}

\section{Proof of the theorems}

From now on, $G$ will be a finite group which has a normal $p$-subgroup 
$N$ containing its centralizer $\cen{G}{N}$.
Still, $R$ is a $p$-adic ring. 
Let $Q$ be a finite $p$-subgroup of $\NU{RG}$ containing $N$ as a normal 
subgroup. We set $M=RG$ and consider $M$ as $R(G\times Q)$-module, the 
action given by $m\cdot r(g,x)=g^{-1}rmx$ for all $g\in G$,
$x\in Q$, $r\in R$ and $m\in M$. We shall write, for example, 
$M_{G\times 1}$ for the restriction of $M$ to $G$ (acting from the left).

We will need the following indecomposability criterion which is proved
in greater generality in \cite[p.~231]{rogg:86}.    
For convenience of the reader, the (short) proof is included.

\begin{lem}\label{ind}
As $R(G\times N)$-module, $M$ is indecomposable.
\end{lem}
\begin{pf}
Any endomorphism of $M_{G\times 1}$ is given by right multiplication with
some element of $RG$, so $\mbox{\rm End}_{R(G\times N)}(M)\cong\cen{RG}{N}$.
The radical quotient $\cen{RG}{N}/\mbox{rad}(\cen{RG}{N})$ is isomorphic to
$\cen{kG}{N}/\mbox{rad}(\cen{kG}{N})$ where $k$ denotes the residue class
field of $R$ (see \cite[(5.22)]{CuRe:81}), and the latter 
quotient is isomorphic to $k$ by Proposition~\ref{sbor}(a) 
since $\AUG{k}{N}G$ is nilpotent (see \cite[(5.26)]{CuRe:81}).
Hence $\mbox{\rm End}_{R(G\times N)}(M)$ 
is local, and $M$ is indecomposable (see \cite[(6.10)]{CuRe:81}).
\qed
\end{pf}

By the Ward--Coleman Lemma, we may fix (and do) for each $x\in Q$ 
some $g_{x}\in G$ such that 
$g_{x}^{-1}x\in\text{\rm C}_{\text{\rm V}(RG)}(N)$.
Set $U=\lara{g_{x}\mid x\in Q}\leq G$. 
Note that for all $x,y\in Q$ we have
$(g_{x}g_{y})^{-1}xy=(g_{y}^{-1}y)(g_{x}^{-1}x)^{y}\in
\text{\rm C}_{\text{\rm V}(RG)}(N)$.
Since $\cen{G}{N}\leq N$, it follows that $x\mapsto g_{x}N$ defines 
a surjective homomorphism $Q\rightarrow UN/N$.
We set $F=\opn{p}{G}$ (the Fitting subgroup of $G$) and $\bar{G}=G/F$.
We shall extend the bar convention when writing $\bar{m}$ for the 
image of $m$ in $RG$ under the natural map $RG\rightarrow R\bar{G}$.

\begin{lem}\label{l2}
Let $x\in Q$. Then $\bar{x}\in \bar{g_{x}}+pR\bar{G}$.
Further, if $g_{x}\in F$ then $\bar{x}=1$. In other words,
$\bar{x}\mapsto\bar{g_{x}}$ defines an isomorphism 
$\bar{Q}\rightarrow\bar{U}$, and $\bar{Q}$ maps onto $\bar{U}$
under the natural map $R\bar{G}\rightarrow (R/pR)\bar{G}$.
\end{lem}
\begin{pf}
Let $x\in Q$. Since $x$ has augmentation one, 
$\bar{x}\in \bar{g_{x}}+pR\bar{G}$ 
by Proposition~\ref{sbor}(a). Now assume that $g_{x}\in F$ and, 
by way of contradiction, that $\bar{x}$ has order $p$.
Then $p=2$ by Lemma~\ref{odd}, and by Proposition~\ref{sbor}(b),
$\bar{x}=\bar{g_{x}}^{-1}\bar{x}\in R+4R\bar{G}+2R[\bar{T}]$
where $\bar{T}$ is the set of involutions in $\bar{G}$.
Hence we may apply Proposition~\ref{cp} 
(to the group $\bar{G}$, with $c=1$ and $u=\bar{x}$) to
conclude that $\bar{x}=1$, a contradiction. The lemma is proved.
\qed
\end{pf}

The next lemma is the place where Weiss' theorem obviously comes into play. 
For a subgroup $H$ of $G$ we shall write $\widehat{H}$ for the sum
of its elements in $RG$. Let $P$ be a Sylow $p$-subgroup of $G$.

\begin{lem}\label{l3}
$M_{P\times Q}$ is a permutation lattice for $P\times Q$ over $R$.
\end{lem}
\begin{pf}
First, we consider the module $\widehat{F}M_{1\times Q}$ which is 
isomorphic to $R\bar{G}$ with $Q$ acting by right multiplication via 
$Q\rightarrow \bar{Q}$. Write $\bar{G}$ as a disjoint union 
$\bar{G}=\bigcup_{i}\bar{h_{i}}\bar{U}$ with $h_{i}\in G$, and set 
$B=\bigcup_{i}\{\bar{h_{i}}\,\bar{x}\mid x\in Q\}$. By Lemma~\ref{l2},
$B$ reduces modulo $pR\bar{G}$ to the canonical basis $\bar{G}$ of
$(R/pR)\bar{G}$, and $|B|=|\bar{G}|$. Thus 
$B$ is an $R$-basis of $R\bar{G}$ by Nakayama's Lemma.
Since the elements of $B$ are permuted under the action of $Q$
it follows that $\widehat{F}M_{1\times Q}$ is a permutation lattice.

Next, we shall show that $\widehat{P}M$ is a permutation lattice for 
$1\times Q$ over $R$. Set $V=\widehat{P}M_{1\times Q}$. Write $G$ as 
a disjoint union $G=\bigcup_{i}Pk_{i}U$ with $k_{i}\in G$ and set
$B=\bigcup_{i}\{\widehat{P}k_{i}x\mid x\in Q\}$.
Write $s$ for the sum of the elements of some system of coset
representatives of $F$ in $P$, so that $\widehat{P}=\widehat{F}s$.
We have $V\subseteq\widehat{F}M\cong R\bar{G}$, and treat the latter 
isomorphism as an identification. So $V=\bar{s}R\bar{G}$ and 
$B=\bigcup_{i}\{\bar{s}\bar{k_{j}}\bar{x}\mid x\in Q\}$.
As above, it follows that $B$ is an $R$-basis of $V$, the 
elements of which are permuted under the action of $Q$.

The lemma thus follows from Weiss' theorem \cite[Theorem~(50.1)]{wei:93}.
\qed
\end{pf}

We are now in a position to prove Theorem~B. We shall make use of the 
elementary theory of vertices and sources, the Krull-Schmidt Theorem and the 
Mackey decomposition. As a general reference we give \cite{CuRe:81}.

\begin{thm}\label{pcg}
If $c$ is a unit in $\NU{RG}$ of finite $p$-power order and $[N,c]=1$, 
then $c\in N$.
\end{thm}
\begin{pf}
Set $C=\lara{c}$, and consider $M$ as $R(G\times N\times C)$-module, the 
action of $G\times N$ being given as before, and 
$m\cdot y=my$ for $m\in M$, $y\in C$.

Note that $C$ acts trivially on $\widehat{N}M$. Indeed, $\widehat{N}M$ is a 
permutation module for $C$ by Lemma~\ref{l3} and $C$ acts trivially on
$\widehat{N}M/p\widehat{N}M$ by Proposition~\ref{sbor}(a).

By Lemma~\ref{l3},
\begin{equation}\label{eq2}
M_{P\times N\times C}\;\cong\;\bigoplus_{j=1}^{n}1\big\uparrow\mbox{}_{U_{j}}
^{P\times N\times C}
\end{equation}
for some subgroups $U_{j}$ of $P\times N\times C$.
The number of summands equals the $R$-rank of the fixed point module
$M^{P\times N\times C}$, which coincides with $M^{P\times 1\times 1}$
($=\widehat{P}RG$). Hence $n=|G:P|$. 
Since the modules $M_{P\times 1\times 1}$ and $M_{1\times N\times 1}$ are
both free, it follows from \eqref{eq2} and Mackey decomposition that
$U_{j}\cap(P\times 1\times 1)=1$
and $U_{j}\cap(1\times N\times 1)=1$ for all $j$.
In particular, $|U_{j}|\leq|N|\cdot|C|$ for all $j$, and comparing ranks in 
\eqref{eq2} gives $|U_{j}|=|N|\cdot|C|$ for all $j$.
Restricting to $N\times N\times C$ in \eqref{eq2} gives
\[ M_{N\times N\times C}\;\cong\;\bigoplus_{j=1}^{|G:P|}\,\bigoplus
_{U_{j}(N\times N\times C)\backslash a} 1\big\uparrow\mbox{}
_{U_{j}^{a}\cap (N\times N\times C)}^{N\times N\times C}. \]
The number of summands in this decomposition is the $R$-rank of
$M^{N\times N\times C}$. 
Since $M^{N\times N\times C}=M^{N\times 1\times 1}$ there are $|G:N|$ 
summands. It follows that $U_{j}\leq N\times N\times C$ for all $j$.
Hence $U_{j}=\{(x,(x\alpha_{j})z_{j}^{-i},c^{i})\mid x\in N, i\in\ZZ\}$
for some $z_{j}\in\zen{N}$ and $\alpha_{j}\in\aut{N}$,
and 
$1\big\uparrow\mbox{}_{U_{j}}^{G\times N\times C}
\big\downarrow\mbox{}_{G\times N\times 1}$ is isomorphic to 
$1\big\uparrow\mbox{}_{U_{j}\cap(G\times N\times 1)}^{G\times N\times 1}$,
that is, to $RG$, the action given by 
$m\cdot(g,x,1)=g^{-1}m(x\alpha_{j}^{-1})$ for all 
$g\in G$, $x\in N$ and $m\in RG$. 
By Lemma~\ref{ind}, it follows that the modules $1\big\uparrow\mbox{}_{U_{j}}
^{G\times N\times C}$ are indecomposable.
Since $M$ is relatively $P\times N\times C$-projective, $M$ is a direct
summand of $M_{P\times N\times C}\big\uparrow\mbox{}^{G\times N\times C}$, 
and by \eqref{eq2}
\[ M_{P\times N\times C}\big\uparrow\mbox{}^{G\times N\times C}
\;\cong\;\bigoplus_{j=1}^{|G:P|}1\big\uparrow\mbox{}_{U_{j}}
^{G\times N\times C}. \]
By Lemma~\ref{ind}, $M$ is indecomposable, so it follows that
$M\cong 1\big\uparrow\mbox{}_{U_{j}}^{G\times N\times C}$ for some $j$. 
Since $(G\times 1\times 1)\cap U_{j}=1$, this shows that there is $u \in M$ 
such that $Gu$ is an $R$-basis of $M$ and
$u=u\cdot(1,z_{j}^{-1},c)=uz_{j}^{-1}c$.
It follows that $u$ is a unit in $RG$ and $c=z_{j}$.
\qed
\end{pf}
\begin{cor}\label{c4}
If $x\in Q$ and $x\in 1+\AUG{R}{N}G$ then $x\in N$.
\end{cor}
\begin{pf}
For such an $x$ we have $g_{x}\in N$ by Proposition~\ref{sbor}(a),
so $x\in N$ follows from Theorem~\ref{pcg}, applied to 
$c=g_{x}^{-1}x\in\cen{Q}{N}$.
\qed
\end{pf}
\begin{lem}\label{l5}
$M_{1\times Q}$ is a free $RQ$-module.
\end{lem}
\begin{pf}
By Lemma~\ref{l3}, $M_{1\times Q}$ is a permutation lattice for $Q$ over $R$.
By Corollary~\ref{c4}, the kernel of the natural map $Q\rightarrow\NU{RG/N}$
is $N$. By Proposition~\ref{sbor}(a), an element $x$ in $Q$ 
acts on $\widehat{N}RG/p\widehat{N}RG$
by right multiplication with $g_{x}$. 
Since $N\leq Q$ and $RG/N$ is free for the multiplication action of 
$UN/N$---recall that $U=\lara{g_{x}\mid x\in Q}$---it follows that
$|Q|=|UN|$ and that each orbit of the action of $Q$ on a basis of 
$M_{1\times Q}$ the elements of which it permutes has length $|Q|$. 
\qed
\end{pf}

Finally we prove Theorem~A, again guided by the remarks from
\cite[p.~231]{rogg:86}.

\begin{thm}\label{new}
There is a unit $u$ of $RG$ with $Q\leq P^{u}$ and $N=N^{u}$.
\end{thm}
\begin{pf}
By Lemma~\ref{l3},
\begin{equation}\label{eq5}
M_{P\times Q}\;\cong\;\bigoplus_{k=1}^{l}1\big\uparrow\mbox{}_{V_{k}}
^{P\times Q} 
\end{equation}
for some subgroups $V_{k}$ of $P\times Q$. 
Since the modules $M_{P\times 1}$ and 
$M_{1\times Q}$ are both free (the latter by Lemma~\ref{l5}), it follows 
from Mackey decomposition that $V_{k}\cap(P\times 1)=1$ and 
$V_{k}\cap(1\times Q)=1$ for all $k$. We shall show that the induced modules 
$1\big\uparrow\mbox{}_{V_{k}}^{G\times Q}$
are indecomposable. Mackey decomposition gives
\[
M_{P\times N}\;\cong\;\bigoplus_{k=1}^{l}\,\bigoplus
_{V_{k}\backslash a/P\times N}
1\big\uparrow\mbox{}_{V_{k}^{a}\cap (P\times N)}^{P\times N}.
\]
If $g_{1},\ldots,g_{|G:P|}$ is a system of right coset representatives of
$P$ in $G$, then
\[ M_{P\times N}\;\cong\;\bigoplus_{j=1}^{|G:P|}1\big\uparrow\mbox{}_{U_{j}}
^{P\times N}\qquad \text{where }
U_{j}=\{(n,n^{g_{j}})\mid n\in N\}. \]
Let $1\leq k\leq l$. From the above, it follows that there is $j=j(k)$ 
such that without lost of generality
$V_{k}\cap(P\times N)=U_{j}$. In particular, $U_{j}\leq V_{k}$.
Altogether, we see that there is a subgroup $Q_{k}$ of $Q$ containing
$N$ so that $V_{k}=\{(x\beta_{k},x)\mid x\in Q_{k}\}$ for some
injective homomorphism $\beta_{k}:Q_{k}\rightarrow P$.
Hence the induced module $1\big\uparrow\mbox{}_{V_{k}}^{G\times Q_{k}}$
is isomorphic to the $R(G\times Q_{k})$-module $RG$, the operation given by 
$m\cdot (g,x)=g^{-1}m(x\beta_{k})$ for all $g\in G$, 
$x\in Q_{k}$ and $m\in RG$, and so is indecomposable by Lemma \ref{ind}
since $N\beta_{k}=N$. Repeated application of Green's Indecomposability 
Theorem \cite[(19.22)]{CuRe:81} now yields that  
$1\big\uparrow\mbox{}_{V_{k}}^{G\times Q}$ is indecomposable. 

Since $M$ is relatively $P\times Q$-projective, 
$M$ is a direct summand of $M_{P\times Q}\big\uparrow\mbox{}^{G\times Q}$,
and by \eqref{eq5}
\[ M_{P\times Q}\big\uparrow\mbox{}^{G\times Q}
\;\cong\;\bigoplus_{k=1}^{l}1\big\uparrow\mbox{}_{V_{k}}^{G\times Q}.\]
By Lemma~\ref{ind}, $M$ is indecomposable. Hence
$M\;\cong\; 1\big\uparrow\mbox{}_{V_{k}}^{G\times Q}$ for some $k$.
Comparing ranks shows that $Q_{k}=Q$.
Since $V_{k}\cap(G\times 1)=1$, this means that there is $u \in M$ such that
$Gu$ is an $R$-basis of $M$, i.e., $u$ is a unit in $RG$, and
$u=u\cdot (x\beta_{k},x)=(x^{-1}\beta_{k})ux$ for
all $x\in Q$. In other words, $Q\leq P^{u}$ and $N=N^{u}$.
\qed
\end{pf}

\bibliographystyle{elsart-num}

\end{document}